\documentclass[12pt,thmsa]{article}
\usepackage{amsmath}
\usepackage{graphicx}
\usepackage{amsfonts}
\usepackage{amssymb}
\usepackage{mathrsfs}

\newtheorem{definition}{Definition}
\newtheorem{remark}{Remark}

\newtheorem{thm}{Theorem}

\begin{document}

\begin{center}
{\Large \textbf{Testing for homogeneity  of  several  functional variables  via multiple  maximum variance discrepancy}}

\bigskip

Armando Sosthène Kali  BALOGOUN$^a$  and  Guy Martial  NKIET$^b$ 

\bigskip

$^a$ Universit\'{e} Dan Dicko Dankoulodo, Maradi, Niger.
 
$^b$ Universit\'{e} des Sciences et Techniques de Masuku,  Franceville, Gabon.

\bigskip

E-mail adresses : armando.balogoun@uddm.edu.ne,    guymartial.nkiet@univ-masuku.org.

\bigskip
\end{center}

\noindent\textbf{Abstract.} This paper adresses  the problem of testing for the equality of $k$ probability distributions on Hilbert spaces, with $k\geqslant 2$. We introduce a generalization of the maximum variance discrepancy called multiple maximum variance discrepancy (MMVD).  Then, a consistent  estimator of  this measure  is proposed as test statistic,  and its asymptotic distribution under the null hypothesis is derived. A simulation study comparing the proposed test with existing ones is provided.
\bigskip


\vskip 4mm

\section{Introduction}\label{Intro}
 Testing for homogenity, that is the equality of several distributions, has long  been a major concern in statistics. Many methods have been proposed in the literature to achieve this in the case where the observations come from samples of real random variables, some  of the most famous being the Komogorov-Smirnov test and  the Cramér-von Mises test. When the data to be processed are more complex such as functional data,  as it is often the case in  number  of applied fields  including environmental science, finance, genetics, biology, geophysics, image and signal processing, it is also important to test for  homogeneity.   However, most of the works considering this problem  in the context of functional data analysis  propose methods which boil down to tests for  the equality of certain characteristics of the distributions  of the concerned functional variables  such as means or covariance operators   (e.g.,  \cite{aghoukeng09,boente18,cuevas04,fremdt12,zhang10}).   On the other hand, the framework of kernel-based methods, that is methods based on kernel embeddings of probability measures in reproducing kernel Hilbert spaces (RKHS), offers the possibility to test for the equality of the distributions of random variables valued into metric spaces, including functional variables. This is how several methods have been proposed in this context, but mainly for the case of two populations. More specifically, the  maximum mean discrepancy (MMD) has been introduced in \cite{gretton12}   for the two-sample problem, that is the test for equality of  distributions of two random variables valued into metric spaces, and a consistent  estimator of it has been used as test statistic. More recently, this has been extended to the case of several populations via the generalized maximum mean discrepancy (GMMD) introduced in \cite{balogoun21,balogoun22}. Further investigations allowed \cite{makigusa23} to propose a novel discrepancy measure, called maximum variance discrepancy (MVD), and to use it for testing the equality of two distributions, so introducing a new testing procedure which is shown, from simulations, to be more powerful than the one based on MMD.    In this paper, we adress the problem of extending the approach of \cite{makigusa23} to the case of more than two populations by introducing the so-called  multiple maximum variance discrepancy (MMVD). More precisely,   for  integers $k\geq 2$ and $j\in\{1,\cdots,k\}$,  considering  a random variable $X_j$    valued into   a separable Hilbert  space $\mathcal{X}$ and with distribution $\mathbb{P}_j$, we are interested in  the problem of testing for the hypothesis
\[
\mathscr{H}_0\,:\,\mathbb{P}_1=\mathbb{P}_2=\cdots=\mathbb{P}_k
\,\,\,\textrm{ 
against }\,\,\,
\mathscr{H}_1\,:\,\,\exists \,(j,\ell)\in\{1,\cdots,k\}^2,\,\,\,\,\,\mathbb{P}_j\neq\mathbb{P}_\ell.
\]
 In Section 2, we define the aforementioned MMVD as an extension of the MVD of \cite{makigusa23}. Then, an estimator of this MMVD is introduced, in Section 3,  as test statistic, and its asymptotic distribution under the null  hypothesis is derived. Section 4, is devoted to the presentation of simulations made in order to evaluate performance of our proposal and to compare it to known  methods. The proof of the main theorem is postponed in Section 5.
 
 %%%%%%%%%%%%%%%%%%%%%%%%%%%%%%%%%%%%%%%%%%%%%%%%%%%%%%%%%%%%%%%%%%%%%%%%%%%%%
 \section{The multiple maximum variance discrepancy}\label{sec2} 
 \noindent In this section,   we introduce the multiple  maximum variance discrepancy (MMVD) as an extension the maximum variance discrepancy (MVD) to the case of more than two populations. 

\bigskip

\noindent Letting $(\mathcal{X}, \mathcal{B})$ be a measurable space, where $\mathcal{X}$ is a separable  Hilbert  space and $\mathcal{B}$ is the corresponding Borel $\sigma$-field, we consider a reproducing kernel Hilbert space (RKHS)  $\mathcal{H}$ of functions from $\mathcal{X}$ to $\mathbb{R}$, endowed with a inner product $<\cdot,\cdot>_{\mathcal{H}}$. That is a Hilbert space for which there exists a map   $K\,:\,\mathcal{X}^2\mapsto\mathbb{R}$, called kernel,  satisfying: $\forall f\in\mathcal{H}$, $f(x)=<K(x,\cdot),f>_\mathcal{H}$.   We assume   that $K$  satisfies the following assumption:

\bigskip
  
 	\noindent $(\mathscr{A}_1)$ :  $\Vert K\Vert_{\infty} := \sup\limits_{(x,y)\in \mathcal{X}^2} K(x,y)<+\infty$. 
\bigskip

 \noindent For $k\geqslant 2$ and $j\in\{1,\cdots, k\}$, we consider a random variable  $X_j$ valued into $\mathcal{X}$  and with probability distribution $\mathbb{P}_j$. Assumption $(\mathscr{A}_1)$ ensures the existence of   the mean and covariance operator embeddings  defined  respectively by:
 \[
 	m_j =\mathbb{E}\left(K\left(X_j, \cdot\right)\right)\,\,\,\textrm{ and }\,\,\,
 	V_j=\mathbb{E}\big(\left(K(X_j,\cdot)-m_j\right)\otimes \left(K(X_j,\cdot)-m_j\right)\big),
 \]
where $\mathbb{E}$ denotes the mathematical expectation  and $\otimes$ is the tensor product such that, for all $(a,b)\in\mathcal{H}^2$, $a\otimes b$ is  the operator from $\mathcal{H}$ to itself  satisfying  $\left(a\otimes b\right)(x)=<a,x>_{\mathcal{H}}b$ for any $x\in \mathcal{H}$. It is known that under  $(\mathscr{A}_1)$, $V_j$ is a Hilbert-Schmidt operator from  $\mathcal{H}$ to itself (see, e.g. \cite{fukumizu07a,fukumizu07b}).  For the case where $k=2$, the maximum variance discrepancy (MVD),    introduced in \cite{makigusa23}  as a measure of the discrepancy between two probability distributions, is defined by:
\begin{equation*}
		\textrm{MVD}(\mathscr{F},\mathbb{P}_1,\mathbb{P}_2):=\sup_{A\in\mathscr{F}}\bigg|<A,V_1-V_2>\bigg|=\left\| V_1-V_2 \right\|,
	\end{equation*}
where $<\cdot,\cdot>$ is the Hilbert-Schmidt  inner product  of operators, $\Vert \cdot\Vert$ is the associated norm and    $\mathscr{F}=\{A\,:\,\parallel A\parallel=1 \}$. We will extend this notion to the case where more than two distributions  are considered, by introducing the multiple maximum variance discrepancy (MMVD):

\bigskip

\begin{definition}\label{d1}
	The multiple maximum variance discrepancy (MMVD), related to $\mathbb{P}_1, \cdots,\mathbb{P}_k$ and $\pi=\left(\pi_1,\cdots,\pi_k\right)\in (]0,1[)^k$ with $\sum_{\ell=1}^k\pi_\ell=1$, is the positive real number defined by			
	\begin{equation*}
		\textrm{MMVD}^2(\mathbb{P}_1, \cdots,\mathbb{P}_k;\pi)=\sum_{j=1}^{k}\sum_{\underset{\ell\neq j}{\ell=1}}^{k}\pi_\ell \,\textrm{MVD}^2(\mathscr{F},\mathbb{P}_j,\mathbb{P}_\ell)=\sum_{j=1}^{k}\sum_{\underset{\ell\neq j}{\ell=1}}^{k}\pi_\ell \parallel V_j-V_\ell\parallel^2\nonumber.
	\end{equation*}
\end{definition} 
\noindent This definition recovers that of  MVD given in \cite{makigusa23} which is clearly a particular case  obtained for $k=2$. If the hypothesis $\mathscr{H}_0$ holds, one has clearly  $\textrm{MMVD}^2(\mathbb{P}_1, \cdots,\mathbb{P}_k;\pi)=0$ for any  $\pi=\left(\pi_1,\cdots,\pi_k\right)\in (]0,1[)^k$ such that $\sum_{\ell=1}^k\pi_\ell=1$. So, a consistent estimator of  $\textrm{MMVD}^2(\mathbb{P}_1, \cdots,\mathbb{P}_k;\pi)$ can be used as test statistic for testing for $\mathscr{H}_0$ against the alternative hypothesis. This is the approach used in \cite{makigusa23}   in the  particular case where $k=2$, and which we take again here in the more general case. 
%%%%%%%%%%%%%%%%%%%%%%%%%%%%%%%%%%%%%%%%%%%%%%%%%%%%%%%%%%%%%%%%%%%%%%%%%%%%%%%%%%%%%%%%%%%%%%%%%%%%%%%%%%%%%%%%%%%%%%%%%%%%%%%%%%%%%%%
\section{Testing for homogeneity}
In this section, the  statistic test  proposed  for testing for $\mathscr{H}_0$ is introduced and we derive its asymptotic distribution under this null  hypothesis.    Finally,   computational aspects are tackled to  show how the introduced test statistic can be computed in practice. 	
 	\subsection{Test statistic} 
 For $j\in\{1,\cdots,k\}$, we consider an i.i.d. sample   $\{X_1^{(j)},\cdots,X_{n_j}^{(j)}\}$ of  $X_{j}$. These samples are such that for $j\neq \ell$,  $X_i^{(j)}$ and $X^{(\ell)}_r$ are independent, for any $(i,r)\in\{1,\cdots,n_j\}\times \{1,\cdots,n_\ell\}$.  In order to estimate the MMVD introduced above, we consider  the empirical estimators of $m_j$ and $V_j$ given by
\begin{eqnarray*}
	\widehat{m}_j=\frac{1}{n_j}\sum_{i=1}^{n_j}K(X_i^{(j)},\cdot)
\end{eqnarray*}
and 
\begin{eqnarray*}
	\widehat{V}_j&=&\frac{1}{n}\sum_{i=1}^{n}\left(K(X_i^{(j)},\cdot)-\widehat{m}_j\right)\otimes \left(K(X_i^{(j)},\cdot)-\widehat{m}_j\right).
\end{eqnarray*} 	
Putting  $n=\sum_{j=1}^{k}n_j$ and  $\pi_j=\frac{n_j}{n}$ for all $j\in \{1,\cdots,k\}$, we assume the following conditions\\
 
\noindent $(\mathscr{A}_2) : $ for any  $j\in \{1,\cdots,k\}$, there exists a real $\rho_j\in]0,1[$ such that   $\lim\limits_{n_j\rightarrow +\infty}\frac{n_j}{n}=\rho_j$ where $\rho_j$.

\bigskip
 
\noindent Then, we take as a statistic test the estimator $\widehat{T}_n$  of
 $\textrm{GMVD}^2(\mathbb{P}_1, \cdots,\mathbb{P}_k;\pi)$ defined by  
 
 \begin{eqnarray*}
 	\widehat{T}_n=\sum_{j=1}^{k}\sum_{\underset{\ell\neq j}{\ell=1}}^{k}\pi_\ell \left\| \widehat{V}_j-\widehat{V}_\ell\right\|^2.
 \end{eqnarray*}  
 
\subsection{Asymptotic distribution under $\mathscr{H}_0$}
 If   $\mathscr{H}_0$ holds, then   $m_1=m_2=\cdots=m_k=:m$ and $V_1=V_2=\cdots=V_k=:V$. Then,  considering  the kernel $\widetilde{K}\, :\,  \mathcal{X}^2\rightarrow \mathbb{R}$ defined by  
 \begin{eqnarray*}
 	\widetilde{K}(x,y)=<\left(K(x,\cdot)-m\right)^{\otimes2}-V,\left(K(y,\cdot)-m\right)^{\otimes2}-V>,
 \end{eqnarray*}
where $a^{\otimes2}=a\otimes a$, and denoting  by $\{\lambda_p\}_{p\geqslant 1}$ the nonincreasing sequence of eigenvalues of the integral operator associated with this kernel, that is the operator   $S_{\widetilde{K}} :L^2(\mathbb{P}) \rightarrow L^2(\mathbb{P})$ defined by
\begin{eqnarray}\label{integralop}
	S_{\widetilde{K}}g(x)=\int_{\mathcal{X}}^{}\widetilde{K}(x,y)g(y)d\mathbb{P}(y),\,\, \textrm{ for }\,\, g\in L^2(\mathbb{P}),
\end{eqnarray}
where $\mathbb{P}=\mathbb{P}_1=\cdots =\mathbb{P}_k$,  we have the following result which gives the asymptotic distribution of the test statistic under null hypothesis:

\bigskip
	
\begin{thm}\label{loi}
	We assume that   the conditions $(\mathscr{A}_1)$  and $(\mathscr{A}_2)$ hold.Then under $\mathscr{H}_0$, as $\min\limits_{1\leqslant j\leqslant k}(n_j)\rightarrow +\infty$,
	\begin{eqnarray*}
		n\widehat{T}_n\stackrel{\mathscr{D}}{\rightarrow}\sum_{p=1}^{+\infty}\lambda_p\left\{(k-2)Z_p+\sum_{j=1}^{k}\left(\rho_j^{-1}Y^2_{j,p}-2\sum_{\underset{\ell\neq j}{\ell=1}}^k\rho_\ell^{1/2}\rho_j^{-1/2}Y_{j,p}Y_{\ell,p}\right)\right\},
	\end{eqnarray*}
where $(Y_{j,p})_{p\geqslant 1,\,1\leqslant j\leqslant k}$ is a sequence of independent random variables having the standard normal distribution  and $(Z_p)_{p\geqslant 1}$ is a sequence of independent random variables having the chi-squared   distribution with $k$ degrees of freedom.
\end{thm} 

\bigskip

\begin{remark}
	This Theorem generalizes  Theorem 3.1 of \cite{makigusa23}. Indeed, if $k=2$ the previous limiting distribution is that of
\[
\frac{1}{\rho (1-\rho)}\sum_{p=1}^{+\infty}\lambda_p\mathcal{W}_p,
\]
where $\rho=\rho_1$ (hence, $\rho_2=1-\rho$), and
\begin{eqnarray*}
\mathcal{W}_p&=&\rho(1-\rho)\Bigg\{\frac{1}{\rho}Y_{1,p}^2+\frac{1}{1-\rho}Y_{2,p}^2-2\bigg(\sqrt{\frac{1-\rho}{\rho}}+\sqrt{\frac{\rho}{1-\rho}}\bigg)Y_{1,p}Y_{2,p}\Bigg\}\\
&=&\rho(1-\rho)\Bigg(\frac{1}{\sqrt{\rho}}Y_{1,p}-  \frac{1}{\sqrt{1-\rho}}Y_{2,p}\Bigg)^2\\
&  &+2\rho(1-\rho)\Bigg\{\frac{1}{\sqrt{\rho(1-\rho)}}-\sqrt{\frac{1-\rho}{\rho}}-\sqrt{\frac{\rho}{1-\rho}}\Bigg\}Y_{1,p}Y_{2,p}\\
&=&\mathcal{Z}^2_p
\end{eqnarray*}
with
\[
\mathcal{Z}_p=\sqrt{\rho(1-\rho)}\Bigg(\frac{1}{\sqrt{\rho}}Y_{1,p}-  \frac{1}{\sqrt{1-\rho}}Y_{2,p}\Bigg).
\]
Since $Y_{1,p}$ and $Y_{2,p}$ are independent with standard normal distribution, it follows that $\mathcal{Z}_p$ has also the standard normal distribution. So, Theorem 3.1 of \cite{makigusa23} is recovered.
\end{remark}

\bigskip

\begin{remark}
	As in \cite{makigusa23}, it is not easy to use the asymptotic distribution obtained in Theorem \ref{loi} for performing the test becuase it is an infinite sum and it is difficult to determine the weights contained in it. So, we will approximate this limiting distribution in Section 4  by using a subsampling method. 
\end{remark}

\subsection{Computational aspects}
In this section, we show how to compute our test  statistic in practice. It is done by using   the kernel trick (see \cite{shawe04}) as it was done in \cite{harchaoui07}. First, notice that our test statistic can be writen as 
\begin{eqnarray}
	\widehat{T}_n
	&=& \sum_{j=1}^{k}\left(1+(k-2)\pi_j\right)  \left\| \widehat{V}_j \right\|^2 -2\sum_{j=1}^{k}\sum_{\underset{\ell\neq j}{\ell=1}}^{k}\pi_\ell<\widehat{V}_j,\widehat{V}_\ell>^2\nonumber
\end{eqnarray}  
For $j=1,\cdots,k$, we consider the operator $G_n^{(j)}$ from $\mathbb{R}^{n_j}$ to $\mathcal{H}$ represented in matrix form as
\begin{align*}
G^{(j)}_{n}=[K(X^{(j)}_{1},.),\cdots,K(X^{(j)}_{n_{j}},.)].
\end{align*}
Then putting $$
G_{n}=[G^{(1)}_{n}G^{(2)}_{n}\cdots G^{(k)}_n],$$ and considering the matrices
\[
 \Lambda_n^{(j,l)} =(G_n^{(j)})^TG_n^{(l)} =\left(
\begin{array}{cccc}
K(X_1^{(j)},X_1^{(l)}) & K(X_1^{(j)},X_2^{(l)}) & \cdots & K(X_1^{(j)},X_{n_l}^{(l)}) \\
K(X_2^{(j)},X_1^{(l)}) & K(X_2^{(j)},X_2^{(l)}) & \cdots &K(X_2^{(j)},X_{n_l}^{(l)}) \\
\vdots & \vdots & \ddots & \vdots \\
K(X_{n_j}^{(j)},X_1^{(l)}) & K(X_{n_j}^{(j)},X_2^{(l)}) & \cdots &K(X_{n_j}^{(j)},X_{n_l}^{(l)}) \\
\end{array}
\right)
\]
and 
\[
\Lambda_n=G_n^TG_n=\left(
\begin{array}{cccc}
\Lambda_n^{(1,1)} & \Lambda_n^{(1,2)} & \cdots & \Lambda_n^{(1,k)} \\
\Lambda_n^{(2,1)} & \Lambda_n^{(2,2)} & \cdots & \Lambda_n^{(2,k)} \\
\vdots & \vdots & \ddots & \vdots \\
\Lambda_n^{(k,1)} & \Lambda_n^{(k,2)} & \cdots & \Lambda_n^{(k,k)} \\
\end{array}
\right),
\]
we put
$$\textrm{\textbf{Q}}_{n_j} =\textrm{\textbf{I}}_{n_j}-n_j^{-1}\textrm{\textbf{1}}_{n_j}\textrm{\textbf{1}}_{n_j}^T,$$
where $\textrm{\textbf{I}}_\ell$ (resp.   $\textrm{\textbf{1}}_\ell$) the $\ell\times \ell$ identity matrix (resp. the $\ell\times 1$ vector whose components are all equal to $1$). Then, as in \cite{harchaoui07}, we can write 
$$\widehat{V}_{j}=n_j^{-1}G_n^{(j)}\textrm{\textbf{Q}}_{n_j}\textrm{\textbf{Q}}_{n_j}^T(G_n^{(j)})^T,$$ 
and, therefore,
\begin{eqnarray*}\label{st}
	\widehat{T}_n&=&\sum_{j=1}^{k}\left(1+(k-2)\pi_j\right)  \left(n_j^{-1}G_n^{(j)}\textrm{\textbf{Q}}_{n_j}\textrm{\textbf{Q}}_{n_j}^T(G_n^{(j)})^T\right)^2\nonumber\\&& -2\sum_{j=1}^{k}\sum_{\underset{\ell\neq j}{\ell=1}}^{k}\pi_\ell\left\{\left(n_j^{-1}G_n^{(j)}\textrm{\textbf{Q}}_{n_j}\textrm{\textbf{Q}}_{n_j}^T(G_n^{(j)})^T\right)\left(n_\ell^{-1}G_n^{(\ell)}\textrm{\textbf{Q}}_{n_\ell}\textrm{\textbf{Q}}_{n_\ell}^T(G_n^{(\ell)})^T\right)^T\right\}\nonumber\\
	&=&\sum_{j=1}^{k}\left(1+(k-2)\pi_j\right)  \left(n_j^{-1}G_n^{(j)}\textrm{\textbf{Q}}_{n_j}\textrm{\textbf{Q}}_{n_j}^T(G_n^{(j)})^T\right)^2\nonumber\\&& -2\sum_{j=1}^{k}\sum_{\underset{\ell\neq j}{\ell=1}}^{k}\frac{\pi_\ell}{n_j n_\ell}\left\{\left(G_n^{(j)}\textrm{\textbf{Q}}_{n_j}\textrm{\textbf{Q}}_{n_j}^T(G_n^{(j)})^T\right)\left(G_n^{(\ell)}\textrm{\textbf{Q}}_{n_\ell}\textrm{\textbf{Q}}_{n_\ell}^T(G_n^{(\ell)})^T\right)\right\}.\nonumber\\
\end{eqnarray*}  
This formula indicates how to practically compute the above introduce test statistic.
 \section{Simulations}

\noindent In this section, we present the results of simulations made in order to evaluate the performances of the homogeneity test based on MMVD, and to compare it to known procedures, namely the method of \cite{balogoun22} based on generalized maximum mean discrepancy and the method of \cite{cuevas04}.  The empirical size and power were computed  through   Monte Carlo simulations, for  the case where $\mathcal{X}=\textrm{L}^2([0,1])$ with three  populations ($k=3$) given by  processes  $X_j(t)$, $j=1,2,3$,   $t\in[0,1]$. Three models were  used for generating the functional data; more specificallly, denoting by $\varepsilon_1(t)$, $\varepsilon_2(t)$, $\varepsilon_3(t)$, $\nu_1(t)$  and $\nu_2(t)$   independent random variables such that $\varepsilon_j(t)\sim  \mathcal{N}(0,t)$, $j=1,2,3$,  $\nu_1(t)\sim  \mathcal{E}xp(t)$ and $\nu_2(t)\sim  \mathcal{P}(t)$ , we considered:

\bigskip

\noindent\textbf{Model 1 : } $X_j(t)=t(1-t)+\varepsilon_j(t)$,  $j=1,2,3$;

\bigskip

\noindent\textbf{Model 2 : } $X_{1}(t)=t(1-t)^5+\varepsilon_1(t)$,\,\,$X_{2}(t)= t^2(1-t)^4+\varepsilon_2(t)$,\,\,$X_{3}(t)= t^3(1-t)^3+ \nu(t)$;

\bigskip

\noindent\textbf{Model 3 : } $X_{1}(t)= t(1-t)^3+\varepsilon_1(t)$,\,\,$X_{2}(t)= t(1-t)^3-t +\nu_2(t)$,\,\,$X_{3}(t)= t(1-t)^3+ \varepsilon_2(t)$.

\bigskip

\noindent  Model 1, which  corresponds to a case where  the  hypothesis $\mathscr{H}_0$ holds,   is   used  for computing empirical size.  For the other two models  $\mathscr{H}_0$ fails; so, they are used for computing empirical power.  We generated  $2000$ independent samples of each of the preceding processes in discretized versions  on equispaced values $t_1,\cdots,t_{21}$  in the interval $[0,1]$, with $t_\ell=(\ell-1)/20$, $\ell=1,\cdots,21$.  The sample sizes where  $n_1=n_2=n_3= 25,\,50,\,100,\,200,\,300$.  For performing  our method and that of \cite{balogoun22}, we used a  Gaussian  kernel $K(x,y)=\exp\left(-0.5\Vert x-y\Vert^2\right)$,    and we  computed the terms
\[
K(X_i^{(j)},X_k^{(\ell)}) =\exp\left(-0.5\int_0^1(X_i^{(j)}(t)-X_k^{(\ell)}(t)) ^2\,\,dt\right)
\]
by approximating the  above integral  by using the trapezoidal rule,  and the $p$-values related to our  test based on MMVD  by the  permutation method using the introduced test statistic.  The nominal significance level was taken as $\alpha=0.05$.   The obtained results are reported in Table 1  in which  our test is   denoted by M1, the one introduced in \cite{balogoun22} is denoted by   M2 and that of  \cite{cuevas04} is  denoted by M3. They  show that all the three methods are well calibrated for finite samples since the obtained empirical sizes are close to the nominal significance level, except for $n=25$ for M3.  Regarding the empirical power, it is seen  that M1 performs very well, and gives comparable results to   M2. These two methods outperform M3 which gives very bad results for Model 3.

\begin{table}
\setlength{\tabcolsep}{0.08cm} %donne la distance entre les collones du tableau%
\renewcommand{\arraystretch}{1} %donne la distance entre les lignes%
\begin{center}
{\begin{tabular}{ccccccccccccccc}
\hline
 &  &  &  &  &  &  &   &  &  &  &  &  &  &     \\
 &  &  &  $n_{1}=n_{2}=n_3 $ &  &  & M1 &  &  & M2 &  &  & M3 &  & \\      
\hline
\hline
 &  &  &  &  &  &  &   &  &  &  &  &  &  &     \\
 &  &  & 25  &  & &  0.060 &  &   & 0.090 &  &  & 0.160  &  &     \\
&  &  & 50  &  & &  0.050 &  &   & 0.060 &  &  & 0.070  &  &     \\
Model 1&  &  & 100  &  & &  0.050 &  &   & 0.040 &  &  & 0.060  &  &     \\
&  &  & 200  &  & &  0.040 &  &   & 0.040 &  &  & 0.050  &  &     \\
&  &  & 300  &  & &  0.050 &  &   & 0.050 &  &  & 0.050  &  &     \\
 &  &  &  &  &  &  &   &  &  &  &  &  &  &     \\    
\hline
 &  &  &  &  &  &  &   &  &  &  &  &  &  &     \\
 &  &  & 25  &  & &  0.960 &  &   & 0.940 &  &  & 0.740  &  &     \\
&  &  & 50  &  & &  0.990 &  &   & 0.990 &  &  & 0.840  &  &     \\
Model 2&  &  & 100  &  & &  1.000 &  &   & 1.000 &  &  & 0.890  &  &     \\
&  &  & 200  &  & &  1.000  &  & &  1.000 &  &  & 0.990  &  &     \\
&  &  & 300  &  & &  1.000  &  & &  1.000 &  &  & 1.000  &  &     \\
 &  &  &  &  &  &  &   &  &  &  &  &  &  &     \\    
\hline
  &  &  &  &  &  &  &   &  &  &  &  &  &  &     \\
 &  &  & 25  &  & &  0.980 &  &   & 0.940 &  &  & 0.210  &  &     \\
&  &  & 50  &  & &  0.990 &  &   & 0.990 &  &  & 0.070  &  &     \\
Model 3&  &  & 100  &  & &  1.000 &  &   & 0.990 &  &  & 0.070  &  &     \\
&  &  & 200  &  & &  1.000  &  & &  1.000 &  &  & 0.060  &  &     \\
&  &  & 300  &  & &  1.000  &  & &  1.000 &  &  & 0.050  &  &     \\
 &  &  &  &  &  &  &   &  &  &  &  &  &  &     \\       
\hline\hline
\end{tabular}}         
\end{center}
\centering \caption{\label{table:tab4}Empirical sizes powers over 2000 replications four our test  and the methods of \cite{balogoun22,cuevas04}, with nominal significance level $\alpha=0.05$}
\end{table}

\section{ Proof of Theorem \ref{loi}}
Putting, for $j\in\{1,\cdots,k\}$,
\begin{eqnarray*}
	\widetilde{V}_j=\frac{1}{n_j}\sum_{i=1}^{n_j}\left(K(X^{(j)}_i,\cdot)-m_j\right)^{\bigotimes2} ,
\end{eqnarray*}
we have $n\widehat{T}_n=A_n + B_n$, where
\begin{eqnarray*}\label{tn}
	A_n
&=&\sum_{j=1}^{k}\sum_{\underset{\ell\neq j}{\ell=1}}^{k}\left(\frac{\pi_\ell}{\pi_j+\pi_\ell}\right)\left\{ (n_j+n_\ell)\left(\left\| \widehat{V}_j-\widehat{V}_\ell\right\|^2-\left\|	\widetilde{V}_j-	\widetilde{V_\ell}\right\|^2\right)\right\},	\\
 B_n&=&n\sum_{j=1}^{k}\sum_{\underset{\ell\neq j}{\ell=1}}^{k}\pi_\ell\left\|	\widetilde{V}_j-	\widetilde{V_\ell}\right\|^2 .\nonumber
\end{eqnarray*}
According to Lemma 1 of \cite{makigusa20}, $(n_j+n_\ell)\left(\left\| \widehat{V}_j-\widehat{V}_\ell\right\|^2-\left\|	\widetilde{V}_j-	\widetilde{V_\ell}\right\|^2\right)=o_P(1)$. Since $\lim\limits_{n_j\rightarrow +\infty}\pi_j=\rho_j$ and $\lim\limits_{n_l\rightarrow +\infty}\pi_\ell=\rho_\ell$, $\lim\limits_{n_l, n_j\rightarrow +\infty}\frac{\pi_\ell}{\pi_j+\pi_\ell}=\frac{\rho_\ell}{\rho_j+\rho_\ell}=\frac{1}{1+\rho_j\rho_\ell^{-1}}$, we deduce that 
$A_n=o_P(1)$ and that  $n\widehat{T}_n$ has the same limiting distribution than $B_n$. It remains to determine this later. We have : 
\begin{eqnarray*}
	& &\left\|	\widetilde{V}_j-	\widetilde{V_\ell}\right\|^2\\
	& &=\frac{1}{n^2_j}\sum_{i=1}^{n_j}\sum_{r=1}^{n_j}<\left(K(X^{(j)}_i,\cdot)-m_j\right)^{\otimes2}-V,\left(K(X^{(j)}_r,\cdot)-m_j\right)^{\otimes2}-V>^2\nonumber\\
	& &+\frac{1}{n^2_\ell}\sum_{i=1}^{n_\ell}\sum_{r=1}^{n_\ell}<\left(K(X^{(\ell)}_i,\cdot)-m_\ell\right)^{\otimes2}-V,\left(K(X^{(\ell)}_r,\cdot)-m_\ell\right)^{\otimes2}-V>^2\nonumber\\
	&-&\frac{2}{n_j n_\ell}\sum_{i=1}^{n_j}\sum_{r=1}^{n_\ell}<\left(K(X^{(j)}_i,\cdot)-m_j\right)^{\otimes2}-V,\left(K(X^{(\ell)}_r,\cdot)-m_\ell\right)^{\otimes2}-V>^2\nonumber\\
&=&
	\frac{1}{n^2_j}\sum_{i=1}^{n_j}\sum_{r=1}^{n_j}\widetilde{K}\left(X_i^{(j)},X_r^{(j)}\right)
	+\frac{1}{n^2_\ell}\sum_{i=1}^{n_\ell}\sum_{r=1}^{n_\ell}\widetilde{K}\left(X_i^{(\ell)},X_r^{(\ell)}\right)\nonumber\\
	&-&\frac{2}{n_j n_\ell}\sum_{i=1}^{n_j}\sum_{r=1}^{n_\ell}\widetilde{K}\left(X_i^{(j)},X_r^{(\ell)}\right).
\end{eqnarray*}
Thus  $B_n=C_n+D_n$, where
\begin{eqnarray}\label{cn}
	C_n&=&\sum_{j=1}^{k} \sum_{\underset{\ell\neq j}{\ell=1}}^{k} \left\{\frac{\pi_j^{-1}\pi_\ell}{n_j}\sum_{i=1}^{n_j}\sum_{r=1}^{n_j}\widetilde{K}\left(X_i^{(j)},X_r^{(j)}\right)\right.\nonumber\\
& &\hspace{3cm}+\left.\frac{1}{n_\ell}\sum_{i=1}^{n_\ell}\sum_{r=1}^{n_\ell}\widetilde{K}\left(X_i^{(\ell)},X_r^{(\ell)}\right)\right\}
\end{eqnarray}
and 
\begin{eqnarray}\label{dn}
	D_n=-2\sum_{j=1}^{k} \sum_{\underset{\ell\neq j}{\ell=1}}^{k}\frac{n\pi_\ell}{n_j n_\ell}\sum_{i=1}^{n_j}\sum_{r=1}^{n_\ell}\widetilde{K}\left(X_i^{(j)},X_r^{(\ell)}\right).
\end{eqnarray}
Since $K$ is bounded, the integral operator  $S_{\widetilde{K}}$  given in \eqref{integralop}  is a Hilbert-Schmidt operator (see Theorem VI.22 in \cite{reed81}). Therefore, from Mercer's theorem (see, e.g., \cite{minh06}), we have:
\begin{eqnarray}\label{deco}
	\widetilde{K}(x,y)=\sum_{p=1}^{+\infty}\lambda_pe_p(x)e_p(y),
\end{eqnarray}
where    $\{e_p\}_{p\geqslant 1}$ is an orthonormal basis of $L^2(\mathbb{P})$ such that $e_p$ is an eigenvector of  $S_{\widetilde{K}}$     associated to   $\lambda_p$. Note that, from  $\int_{\chi}^{}\widetilde{K}(x,y)e_p(y)d\mathbb{P}(y)=\lambda_pe_p(x)$, we get:
\begin{eqnarray*}
	\mathbb{E}\left(e_p(X^{(j)}_i)\right)=\frac{1}{\lambda_p}\int_{\chi}^{}\mathbb{E}\left(\widetilde{K}(X^{(j)}_i,y)\right)e_p(y)d\mathbb{P}(y),
\end{eqnarray*}
and since
\begin{eqnarray*}
	&&\mathbb{E}\left(\widetilde{K}(X^{(j)}_i,y)\right)\\
&&=<\mathbb{E}\left(\left(K(X^{(\ell)},\cdot)-m_\ell\right)^{\otimes2}-V\right),\left(K(y,\cdot)-m_\ell\right)^{\otimes2}-V>=0
\end{eqnarray*} 
it follows that $\mathbb{E}\left(e_p(X^{(j)}_i)\right)=0$.  From \eqref{cn} and \eqref{deco},
\begin{eqnarray*}
	C_n
	&=&\sum_{j=1}^{k} \sum_{\underset{\ell\neq j}{\ell=1}}^{k} \left\{\frac{\pi_j^{-1}\pi_\ell}{n_j}\sum_{i=1}^{n_j}\sum_{r=1}^{n_j}\sum_{p=1}^{+\infty}\lambda_pe_p(X_i^{(j)})e_p(X_r^{(j)})\right\}\nonumber\\
& &+\sum_{j=1}^{k} \sum_{\underset{\ell\neq j}{\ell=1}}^{k} \left\{ \frac{1}{n_\ell}\sum_{i=1}^{n_\ell}\sum_{r=1}^{n_\ell}\sum_{p=1}^{+\infty}\lambda_pe_p(X_i^{(\ell)})e_p(X_r^{(\ell)})\right\}\nonumber\\
&=&\sum_{p=1}^{+\infty}\lambda_p\left\{\sum_{j=1}^{k}  \pi_j^{-1}\left(\sum_{\underset{\ell\neq j}{\ell=1}}^{k}\pi_\ell\right)U_{n_j,p}^2\right\}+\sum_{j=1}^{k} \sum_{\underset{\ell\neq j}{\ell=1}}^{k}  \sum_{p=1}^{+\infty}\lambda_pU_{n_\ell,p}^2,\nonumber\\
\end{eqnarray*}
where
\begin{equation*}\label{unjp}
U_{n_j,p}=\frac{1}{\sqrt{n_j}}\sum_{i=1}^{n_j}e_p(X_i^{(j)}).
\end{equation*}
Since $\sum_{\underset{\ell\neq j}{\ell=1}}^{k}\pi_\ell=1-\pi_j$, it follows
\begin{eqnarray}\label{cnfinal}
	C_n&=&\sum_{p=1}^{+\infty}\lambda_p\left\{\sum_{j=1}^{k}  \pi_j^{-1}\left(1-\pi_j\right)U_{n_j,p}^2\right\}+\sum_{p=1}^{+\infty}\lambda_p\left\{\sum_{j=1}^{k} \sum_{\underset{\ell\neq j}{\ell=1}}^{k}  U_{n_\ell,p}^2\right\}\nonumber\\
	&=&\sum_{p=1}^{+\infty}\lambda_p\left\{\sum_{j=1}^{k}  \left(\pi_j^{-1}+k-2\right)U_{n_j,p}^2\right\}.
\end{eqnarray}
Further, from \eqref{dn} and \eqref{deco},
\begin{eqnarray*}
	D_n	
	&=&-2\sum_{p=1}^{+\infty}\lambda_p\sum_{j=1}^{k}\frac{\pi_j^{-1/2}}{\sqrt{n_j}}\sum_{i=1}^{n_j}e_{p}(X_i^{(j)})\left(\sum_{\underset{\ell\neq j}{\ell=1}}^k\frac{\pi_\ell^{1/2}}{\sqrt{n_\ell}}\sum_{r=1}^{n_\ell}e_{p}(X_r^{(\ell)})\right)\nonumber\\
&=&-2\sum_{p=1}^{+\infty}\lambda_p\sum_{j=1}^{k}\left\{\frac{\pi_j^{-1/2}}{\sqrt{n_j}}\sum_{i=1}^{n_j}e_{p}(X_i^{(j)})\left(\sum_{\ell=1}^{k}\frac{\pi_\ell^{1/2}}{\sqrt{n_\ell}}\sum_{r=1}^{n_\ell}e_{p}(X_r^{(\ell)}) \right.\right.\nonumber\\
&&\left.\left.\hspace{7cm}-\frac{\pi_j^{1/2}}{\sqrt{n_j}}\sum_{r=1}^{n_j}e_p(X_r^{(j)})\right)\right\}\nonumber\\
	&=& -2\sum_{p=1}^{+\infty}\lambda_p\left(\sum_{j=1}^{k}\pi_j^{-1/2}U_{n_j,p}\right)\left(\sum_{\ell=1}^{k}\pi_\ell^{1/2}U_{n_\ell,p}\right)+2\sum_{p=1}^{\infty}\lambda_p\sum_{j=1}^{k}U_{n_j,p}^2.
\end{eqnarray*}
Then, using the equalities $-2ab=(a-b)^2-a^2-b^2$ and $\pi_j^{-1/2}-\pi_j^{1/2}=\frac{1-\pi_j}{\sqrt{\pi_j}}$, it follows
\begin{eqnarray}\label{dnfinal}
	D_n&=&\sum_{p=1}^{+\infty}\lambda_p\left\{2\sum_{j=1}^{k}U_{n_j,p}^2+\left(\sum_{j=1}^{k}\frac{1-\pi_j}{\sqrt{\pi_j}}U_{n_j,p}\right)^2\right\}\nonumber\\
	&&-\sum_{p=1}^{+\infty}\lambda_p\left\{\left(\sum_{j=1}^{k}\pi_j^{-1/2}U_{n_j,p}\right)^2+\left(\sum_{j=1}^{k}\pi_j^{1/2}U_{n_j,p}\right)^2\right\}.
\end{eqnarray}
Using \eqref{cnfinal} and \eqref{dnfinal}, we get
\begin{eqnarray*}
	B_n
	&=&\sum_{p=1}^{+\infty}\lambda_p\left\{\sum_{j=1}^{k}  \left(\pi_j^{-1}+k\right)U^2_{n_j,p}+\left(\sum_{j=1}^{k}(\pi_j^{-1/2}-\pi_j^{1/2})U_{n_j,p}\right)^2\right\}\nonumber\\
	&-&\sum_{p=1}^{+\infty}\lambda_p\left\{\left(\sum_{j=1}^{k}\pi_j^{-1/2}U_{n_j,p}\right)^2+\left(\sum_{j=1}^{k}\pi_j^{1/2}U_{n_j,p}\right)^2\right\}\nonumber\\
	&=&\sum_{p=1}^{+\infty}\lambda_p\Phi_n(\mathcal{U}_{n,p}),
\end{eqnarray*}
where $\mathcal{U}_{n,p}=\left(U_{n_1,p},\cdots,U_{n_k,p}\right)$
and $\Phi_n\,:\,\mathbb{R}^p\rightarrow \mathbb{R}$ is the map defined by
\begin{eqnarray*}
	\Phi_n(x)&=&\sum_{j=1}^{k}  \left(\pi_j^{-1}+k\right)x^2_j+\left(\sum_{j=1}^{k}(\pi_j^{-1/2}-\pi_j^{1/2})x_j\right)^2\\
&&-\left(\sum_{j=1}^{k}\pi_j^{-1/2}x_j\right)^2-\left(\sum_{j=1}^{k}\pi_j^{1/2}x_j\right)^2,
\end{eqnarray*}
where $x_j$ is the $j$-th component of  $x\in\mathbb{R}^p$. Since for any $(j,\ell)\in\{1,2,\cdots,k\}^2$ with  $j\neq\ell$, $U_{n_j,p}$ and   $U_{n_\ell,p}$ are independent,  we get by  the central limit theorem,   $U_{n_j,p}\stackrel{\mathscr{D}}{\rightarrow}Y_{j,p}$ as $n_j\rightarrow +\infty$ where  $Y_{j,p}\sim \mathcal{N}(0,1)$, and $Y_{j,p}$ and $Y_{\ell,p}$ are independent if $j\neq \ell$. Consequently,    $\mathcal{U}_{n,p}\stackrel{\mathscr{D}}{\rightarrow}\mathcal{U}_{p}:=(Y_{1,p},\cdots,Y_{k,p})$ as $\min\limits_{1\leqslant j\leqslant k}(n_j)\rightarrow +\infty$. So, considering the map $\Phi\,:\,\mathbb{R}^p\rightarrow \mathbb{R}$   defined by
\begin{eqnarray*}
	\Phi(x)&=&\sum_{j=1}^{k}  \left(\rho_j^{-1}+k\right)x^2_j+\left(\sum_{j=1}^{k}(\rho_j^{-1/2}-\rho_j^{1/2})x_j\right)^2\\
&&-\left(\sum_{j=1}^{k}\rho_j^{-1/2}x_j\right)^2-\left(\sum_{j=1}^{k}\rho_j^{1/2}x_j\right)^2,
\end{eqnarray*}
we will show that  
\begin{equation}\label{aprouver}
B_n \stackrel{\mathscr{D}}{\rightarrow}\sum\limits_{p=1}^{+\infty}\lambda_p\,\Phi(\mathcal{U}_p)\,\,\,\textrm{ as }\,\,\,
\min_{1\leqslant j\leqslant k}(n_j)\rightarrow +\infty,
\end{equation}
what leads to the required result since
\begin{eqnarray*}
	\Phi(\mathcal{U}_p)&=&\sum_{j=1}^{k}\left(\rho_j^{-1}+k\right)Y_{j,p}^2+\left(\sum_{j=1}^{k}(\rho_j^{-1/2} -\rho_j^{1/2})Y_{j,p}\right)^2\\
&&-\left(\sum_{j=1}^{k}\rho_j^{-1/2}Y_{j,p}\right)^2-\left(\sum_{j=1}^{k}\rho_j^{1/2}Y_{j,p}\right)^2\nonumber\\
	&=&\sum_{j=1}^{k}\left(\rho_j^{-1}+k-2\right)Y_{j,p}^2-2\sum_{j=1}^{k}\sum_{\underset{\ell\neq j}{\ell=1}}^k\rho_\ell^{1/2}\rho_j^{-1/2}Y_{j,p}Y_{\ell,p}\nonumber\\
	&=&(k-2)Z_p+\sum_{j=1}^{k}\left(\rho_j^{-1}Y^2_{j,p}-2\sum_{\underset{\ell\neq j}{\ell=1}}^k\rho_\ell^{1/2}\rho_j^{-1/2}Y_{j,p}Y_{\ell,p}\right),\nonumber
\end{eqnarray*}
where  $Z_{p}=\sum_{j=1}^kY_{j,p}^2\sim\mathcal{X}^2_k$. 

\bigskip

\noindent We denote by $\varphi_U$ the characteristic function of the random variable $U$. The proof of \eqref{aprouver} will be obtained from three steps.

\bigskip 

\noindent\textbf{First step:} We put  $S_n^{(q)}=\sum\limits_{p=1}^{q}\lambda_p\,\Phi(\mathcal{U}_{n,p})$ and we show that for all $\varepsilon>0$ and all $t\in\mathbb{R}$,  there exists $q\in \mathbb{N}^{*}$ such that $\left\vert\varphi_{B_n}(t)-\varphi_{S_n^{(q)}}(t)\right\vert<\frac{\varepsilon}{3}$ for $n_j$ large enough ($j=1,\cdots,k$).

\bigskip

\noindent Using the inequality $\vert e^{iz}-1\vert\leqslant\vert z\vert$ for all $z\in\mathbb{R}$, we have for any $t\in\mathbb{R}$:
\begin{eqnarray*}
	\left\vert\varphi_{B_n}(t)-\varphi_{S_n^{(q)}}(t)\right\vert
&\leqslant& \mathbb{E}\left(\left\vert e^{itB_n}-e^{itS_n^{(q)}}\right\vert\right)\leqslant\mathbb{E}\left(\left\vert e^{it\left(B_n-S_n^{(q)}\right)}-1\right\vert\right)\\
	&\leqslant&\left\vert t\right\vert\mathbb{E}\left(\left\vert  B_n-S_n^{(q)} \right\vert\right)\leqslant\left\vert t\right\vert\sum\limits_{p=q+1}^{+\infty}\lambda_p\mathbb{E}\left(\left\vert  \Phi_n(\mathcal{U}_{n,p}) \right\vert\right).\nonumber
\end{eqnarray*}
On the other hand, we have
\begin{eqnarray*}
	\mathbb{E}\left(\left\vert  \Phi_n(\mathcal{U}_{n,p}) \right\vert\right)
	&\leqslant&\sum_{j=1}^{k}\left(\pi_j^{-1}+k\right)\mathbb{E}\left(U_{n_j,p}^2\right)+\mathbb{E}\left(\left(\sum_{j=1}^{k}(\pi_i^{-1/2}-\pi_j^{1/2})U_{n_j,p}\right)^2\right)\nonumber\\
	& &+\mathbb{E}\left(\left(\sum_{j=1}^{k}\pi_j^{-1/2}U_{n_j,p}\right)^2\right)+\mathbb{E}\left(\left(\sum_{j=1}^{k}\pi_j^{1/2}U_{n_j,p}\right)^2\right).\nonumber
\end{eqnarray*}
Since  $\mathbb{E}\left(e_p^2(X_i^{(j)})\right)=1$  and $\mathbb{E}\left(e_p(X_i^{(j)})\,e_p(X_r^{(\ell)})\right)=\delta_{ir}\delta_{j\ell}$, it follows:
\begin{eqnarray}
	\mathbb{E}\left(U_{n_j,p}^2\right)=\frac{1}{n_j}\sum_{i=1}^{n_j}\mathbb{E}\left(e_p^2(X_i^{(j)})\right)+\frac{1}{n_j}\sum_{i=1}^{n_j}\sum_{\underset{r\neq i}{r=1}}^{n_j}
	\mathbb{E}\left(e_p(X_i^{(j)})\,e_p(X_r^{(j)})\right)=1\nonumber
\end{eqnarray} 
and for $j\neq\ell$,
\[
\mathbb{E}\left(U_{n_j,p}U_{n_\ell,p}\right)=\frac{1}{\sqrt{n_jn_\ell}}\sum_{i=1}^{n_j}\sum_{r=1}^{n_\ell}
	\mathbb{E}\left(e_p(X_i^{(j)})\,e_p(X_r^{(\ell)})\right)=0.
\]
Therefore,
\begin{eqnarray}
	\mathbb{E}\left(\left(\sum_{j=1}^{k}\left(\pi_j^{-1/2}-\pi_j^{1/2}\right)U_{n_j,p}\right)^2\right)&=&\sum_{j=1}^{k}\frac{(1-\pi_j)^2}{\pi_j}\mathbb{E}\left(U_{n_j,p}^2\right)\nonumber\\
	& &+\sum_{j=1}^{k}\sum_{\underset{\ell\neq j}{\ell=1}}^{k}\frac{(1-\pi_j)(1-\pi_\ell)}{\sqrt{\pi_j\pi_\ell}}\mathbb{E}\left(U_{n_j,p}U_{n_\ell,p}\right)\nonumber\\
	&=&\sum_{j=1}^{k}\frac{(1-\pi_j)^2}{\pi_j},\nonumber
\end{eqnarray} 
\begin{eqnarray}
	\mathbb{E}\left(\left(\sum_{j=1}^{k}\pi_j^{-1/2}U_{n_j,p}\right)^2\right)&=&\sum_{j=1}^{k} \pi_j^{-1} \mathbb{E}\left(U_{n_j,p}^2\right) +\sum_{j=1}^{k}\sum_{\underset{\ell\neq j}{\ell=1}}^{k} \pi_j^{-1} \pi_\ell^{-1}\mathbb{E}\left(U_{n_j,p}U_{n_\ell,p}\right)\nonumber\\
	&=&\sum_{j=1}^{k} \pi_j^{-1},\nonumber
\end{eqnarray}
\begin{eqnarray}
	\mathbb{E}\left(\left(\sum_{j=1}^{k}\pi_j^{1/2}U_{n_j,p}\right)^2\right)&=&\sum_{j=1}^{k} \pi_j\mathbb{E}\left(U_{n_j,p}^2\right) +\sum_{j=1}^{k}\sum_{\underset{\ell\neq j}{\ell=1}}^{k} \pi_j^{1/2} \pi_\ell^{1/2}\mathbb{E}\left(U_{n_j,p}U_{n_\ell,p}\right)\nonumber\\
	&=&\sum_{j=1}^{k} \pi_j=1,\nonumber
\end{eqnarray}
and, consequently,
\begin{eqnarray*}
	\mathbb{E}\left(\left\vert  \Phi_n(\mathcal{U}_{n,p}) \right\vert\right)
	&\leqslant&\sum_{j=1}^{k}\left(\pi_j^{-1}+k+\frac{(1-\pi_j)^2+1}{\pi_j}+\pi_j\right)\nonumber\\
	&=&\sum_{j=1}^{k}\left(3\pi_j^{-1}+2\pi_j+k-2\right).\nonumber\\
\end{eqnarray*} 
Since $\lim_{n_j\rightarrow +\infty}\left(3\pi_j^{-1}+2\pi_j+k-2\right)=3\rho_j^{-1}+2\rho_j+k-2$, there exists $n_j^0\in\mathbb{N}^\ast$ such that, for any $n_j\geqslant n_j^0$, one has $3\pi_j^{-1}+2\pi_j+k-2\leqslant 3\rho_j^{-1}+2\rho_j+k-1$. Hence, for $n_1\geqslant n_1^0,\cdots,n_k\geqslant n_k^0$,
\begin{eqnarray}
	\mathbb{E}\left(\left\vert  \Phi_n(\mathcal{U}_{n,p}) \right\vert\right)
	&\leqslant&\sum_{j=1}^{k}\left(3\rho_j^{-1}+2\rho_j+k-1\right),\nonumber
\end{eqnarray}
and, consequently,
\begin{eqnarray}
	\left\vert\varphi_{B_n}(t)-\varphi_{S_n^{(q)}}(t)\right\vert
	\leqslant\left\vert t\right\vert\sum_{j=1}^{k}\left(3\rho_j^{-1}+2\rho_j+k-2\right)\sum\limits_{p=q+1}^{+\infty}\lambda_p.\nonumber
\end{eqnarray}
Since $\sum_{p=1}^{+\infty}\lambda_p<+\infty$, one has $\lim_{q\rightarrow +\infty}\left\vert\varphi_{B_n}(t)-\varphi_{S_n^{(q)}}(t)\right\vert=0$. So, there exists $q_0\in\mathbb{N}^\ast$ such that, for any $q\geqslant q_0$, 
\begin{equation}\label{eps1}
	\left\vert\varphi_{B_n}(t)-\varphi_{S_n^{(q)}}(t)\right\vert<\frac{\varepsilon}{3}.
\end{equation}

\bigskip

\noindent\textbf{Second step:} We consider $S_q=\sum_{p=1}^q\lambda_p\Phi(\mathcal{U}_p)$ and we show that we have   $S_n^{(q)}\stackrel{\mathscr{D}}{\rightarrow}S_q$  as $\min\limits_{1\leqslant j\leqslant k}(n_j)\rightarrow +\infty$.\\

\noindent It suffices to show that $S_n^{(q)}-S_q$ converges in probability to $0$. We have the following inequality :

\begin{eqnarray}\label{ineq1}
	\left\vert S_n^{(q)}-S_q\right\vert
	&\leqslant&\sum_{p=1}^q\lambda_p\left\vert\Phi_n(\mathcal{U}_{n_j,p})-\Phi(\mathcal{U}_{p})\right\vert\nonumber\\
	&\leqslant &\sum_{p=1}^q\lambda_p\bigg(\left\vert\Phi_n(\mathcal{U}_{n,p})-\Phi(\mathcal{U}_{n,p})\right\vert+\left\vert\Phi(\mathcal{U}_{n,p})-\Phi(\mathcal{U}_p)\right\vert\bigg).
\end{eqnarray}
Moreover, by using $a^2-b^2=(a-b)^2+2b(a-b)$, we get
\begin{eqnarray}\label{ineq3}
	&&\left\vert\Phi_n(\mathcal{U}_{n,p})-\Phi(\mathcal{U}_{n,p})\right\vert\nonumber\\
	& &\leqslant \bigg\{\sum_{j=1}^{k}\left\vert\pi_j^{-1}-\rho^{-1}_j\right\vert\nonumber\\
	&&+\left(\sum_{j=1}^{k}\left\vert\frac{1-\pi_j}{\sqrt{\pi_j}}-\frac{1-\rho_j}{\sqrt{\rho_j}}\right\vert\right)^2
	+ 2\sum_{j=1}^{k} \sum_{\ell=1}^{k}\frac{1-\rho_\ell}{\sqrt{\rho_\ell}}\left\vert\frac{1-\pi_j}{\sqrt{\pi_j}}-\frac{1-\rho_j}{\sqrt{\rho_j}}\right\vert\nonumber\\
	& &+\left(\sum_{j=1}^{k}\left\vert \pi_j^{-1/2}- \rho_j^{-1/2} \right\vert\right)^2
	+ 2\sum_{j=1}^{k} \sum_{\ell=1}^{k} \rho_\ell^{-1/2}\left\vert \pi_j^{-1/2}- \rho_j^{-1/2}\right\vert\nonumber\\
	& &+\left(\sum_{j=1}^{k}\left\vert \pi_j^{1/2}- \rho_j^{1/2} \right\vert\right)^2
	+ 2\sum_{j=1}^{k} \sum_{\ell=1}^{k} \rho_\ell^{1/2}\left\vert \pi_j^{1/2}- \rho_j^{1/2}\right\vert\,\bigg\}\left\Vert \mathcal{U}_{n,p}\right\Vert^2_{\mathbb{R}^p},
\end{eqnarray}
where $\Vert\cdot\Vert_{\mathbb{R}^p}$ denotes the Euclidean norm of $\mathbb{R}^p$. Since $\mathcal{U}_{n,p}$ converges in distribution  to $\mathcal{U}_{p}$, we deduce from (\ref{ineq1}), (\ref{ineq3}) and of the continuity of $\Phi$ that  $S_n^{(q)}-S_q$ converges in probability to $0$ as $\min\limits_{1\leqslant j\leqslant k}(n_j)\rightarrow +\infty$ and, therefore, that $S_n^{(q)}\stackrel{\mathscr{D}}{\rightarrow}S_q $ as $\min\limits_{1\leqslant j\leqslant k}(n_j)\rightarrow +\infty $.  Thus, there exists $N_1$ such that, for
$\min\limits_{1\leqslant j\leqslant k}(n_j)\geqslant N_1$, we have
\begin{equation}\label{eps2}
	\left\vert\varphi_{S_n^{(q)}}(t)-\varphi_{S_{q}}(t)\right\vert<\frac{\varepsilon}{3}.
\end{equation}

\noindent\textbf{Third step:} We put  $S=\sum_{p=1}^{+\infty}\lambda_p\Phi(\mathcal{U}_{p})$ and we show that    $S_{q}\stackrel{\mathscr{D}}{\rightarrow}S$ as  $q\rightarrow +\infty$.\\

\noindent We have:
\begin{eqnarray}
	\left\vert\varphi_{S_q}(t)-\varphi_{S}(t)\right\vert&\leqslant& 
	\mathbb{E}\left(\left\vert e^{itS_q}-e^{itS}\right\vert\right)\leqslant\left\vert t\right\vert\mathbb{E}\left(\left\vert  S_q-S \right\vert\right)
	\leqslant\left\vert t\right\vert\sum\limits_{p=q+1}^{+\infty}\lambda_p\mathbb{E}\left(\left\vert  \Phi(\mathcal{U}_p) \right\vert\right)\nonumber
\end{eqnarray}
and
\begin{eqnarray}
	\mathbb{E}\left(\left\vert  \Phi(\mathcal{U}_p) \right\vert\right)
	&\leqslant&\sum_{j=1}^{k}\left(\rho_j^{-1}+k\right)\,\mathbb{E}\left(Y_{j,p}^2\right)+\mathbb{E}\left(\left(\sum_{j=1}^{k}\frac{1-\rho_j}{\sqrt{\rho_j}}Y_{j,p}\right)^2\right)\nonumber\\
	& &+\mathbb{E}\left(\left(\sum_{j=1}^{k}\rho_j^{-1/2}Y_{j,p}\right)^2\right)+\mathbb{E}\left(\left(\sum_{j=1}^{k}\rho_j^{1/2}Y_{j,p}\right)^2\right).\nonumber
\end{eqnarray}
Since  $\mathbb{E}\left(Y_{j,p}Y_{\ell,p}\right)=\delta_{j\ell}$, it follows
\begin{eqnarray}
	\mathbb{E}\left(\left(\sum_{j=1}^{k}\frac{1-\rho_j}{\sqrt{\rho_j}}Y_{j,p}\right)^2\right)&=&\sum_{j=1}^{k}\frac{(1-\rho_j)^2}{\rho_j}\mathbb{E}\left(Y_{j,p}^2\right)\nonumber\\
	& &+\sum_{j=1}^{k}\sum_{\underset{\ell\neq j}{\ell=1}}^{k}\frac{(1-\rho_j)(1-\rho_\ell)}{\sqrt{\rho_j\rho_\ell}}\mathbb{E}\left(Y_{j,p}Y_{\ell,p}\right)\nonumber\\
	&=&\sum_{j=1}^{k}\frac{(1-\rho_j)^2}{\rho_j},\nonumber
\end{eqnarray} 
\begin{eqnarray}
	\mathbb{E}\left(\left(\sum_{j=1}^{k}\rho_j^{-1/2}Y_{j,p}\right)^2\right)&=&\sum_{j=1}^{k} \rho_j^{-1} \mathbb{E}\left(Y_{j,p}^2\right) +\sum_{j=1}^{k}\sum_{\underset{\ell\neq j}{\ell=1}}^{k} \rho_j^{-1/2} \rho_\ell^{-1/2}\mathbb{E}\left(Y_{j,p}Y_{\ell,p}\right)\nonumber\\
	&=&\sum_{j=1}^{k} \rho_j^{-1} ,\nonumber
\end{eqnarray}
and
\begin{eqnarray}
	\mathbb{E}\left(\left(\sum_{j=1}^{k}\rho_j^{1/2}Y_{j,p}\right)^2\right)&=&\sum_{j=1}^{k} \rho_j \mathbb{E}\left(Y_{j,p}^2\right) +\sum_{j=1}^{k}\sum_{\underset{\ell\neq j}{\ell=1}}^{k} \rho_j^{1/2} \rho_\ell^{1/2}\mathbb{E}\left(Y_{j,p}Y_{\ell,p}\right)\nonumber\\
	&=&\sum_{j=1}^{k} \rho_j=1 .\nonumber
\end{eqnarray}
Hence
\begin{eqnarray}
	\mathbb{E}\left(\left\vert  \Phi(\mathcal{U}_p) \right\vert\right)
	&\leqslant&\sum_{j=1}^{k}\bigg\{3\rho_j^{-1}+2\rho_j+k-2\bigg\}\nonumber
\end{eqnarray}
and
\begin{eqnarray}
	\left\vert\varphi_{S_q}(t)-\varphi_{S}(t)\right\vert&\leqslant& 
	\sum_{j=1}^{k}\bigg\{3\rho_j^{-1}+2\rho_j+k-2\bigg\}\sum\limits_{p=q+1}^{+\infty}\lambda_p.\nonumber
\end{eqnarray}
Since $\sum_{p=1}^{+\infty}\lambda_p<+\infty$, we deduce that $\lim_{q\rightarrow +\infty}\left\vert\varphi_{S_q}(t)-\varphi_{S}(t)\right\vert=0$. So there exists $q_1\in\mathbb{N}^\ast$ such that, for any $q\geqslant q_1$,
\begin{equation}\label{eps3}
	\left\vert\varphi_{S_{q}}(t)-\varphi_{S}(t)\right\vert<\frac{\varepsilon}{3}.
\end{equation}

\bigskip

\noindent Finally, putting $u=\max(q_0,q_1)$, $N_0=\max(n_1^0,\cdots,n_k^0)$  and using (\ref{eps1}), (\ref{eps2}) and (\ref{eps3}), we deduce that if  $\min_{1\leqslant j\leqslant k}(n_j)\geqslant \max(N_0,N_1)$ then
\begin{eqnarray}
	\left\vert\varphi_{B_n}(t)-\varphi_{S}(t)\right\vert\leqslant \left\vert\varphi_{B_n}(t)-\varphi_{S_n^{(u)}}(t)\right\vert
	+\left\vert\varphi_{S_n^{(u)}}(t)-\varphi_{S_{u}}(t)\right\vert
	+\left\vert\varphi_{S_{u}}(t)-\varphi_{S}(t)\right\vert <\varepsilon.\nonumber
\end{eqnarray}
This shows that \eqref{aprouver} holds.

\end{document}